\documentclass[11pt, a4paper]{article}

\usepackage{amssymb}
\usepackage{amsmath}
\usepackage{amscd}
\usepackage{theorem}
\usepackage{mathrsfs}

\usepackage{graphicx}
\usepackage[T1]{fontenc}
\usepackage{ae}
\usepackage{pictexwd,dcpic}

\theoremstyle{change}

\newtheorem{thm}{Theorem.}
\newtheorem{lem}[thm]{Lemma.}

\newtheorem{cor}[thm]{Corollary.}

\renewcommand{\phi}{\varphi}

\newcommand{\proof}{\textbf{Proof. }}
\newcommand{\proofend}{\hfill $\Box$}

\begin{document}

\centerline{\textbf{\large Bergman inner functions and $m$-hypercontractions}}


\vspace{0.2cm}

\centerline{J\"org Eschmeier}
\vspace{1cm}

\begin{center}
\parbox{11cm}{\small
Let $H_m(\mathbb B,\mathcal D)$ be the $\mathcal D$-valued functional Hilbert space
with reproducing kernel $K_m(z,w) = (1-\langle z,w\rangle)^{-m}1_{\mathcal D}$. 
A $K_m$-inner function is by definition an
operator-valued analytic function $W: \mathbb B \rightarrow L(\mathcal E, \mathcal D)$ such that
$\|Wx\|_{H_m(\mathbb B,\mathcal D)} = \|x\|$ for all $x \in \mathcal E$ and $(W\mathcal E) \perp M_z^{\alpha}(W\mathcal E)$
for all $\alpha \in \mathbb N^n \setminus \{0\}$. We show that 
the $K_m$-inner functions are precisely 
the functions of the form $W(z) = D + C \sum^m_{k=1}(1 - ZT^*)^{-k}ZB$,
where $T \in L(H)^n$ is a pure $m$-hypercontraction and the operators $T^*, B, C,D$ form a $2 \times 2$-operator
matrix satisfying suitable conditions. Thus we extend results proved by Olofsson on the unit disc to the case of the unit ball $\mathbb B \subset \mathbb C^n$.
}
\end{center}
\vspace{0.5cm}

\centerline{\textbf{\S1\, Introduction}}

A commuting tuple $T = (T_1, \ldots , T_n) \in L(H)^n$ of bounded linear operators on a complex Hilbert space $H$ is by 
definition a row contraction if the operator
\[
H^n \rightarrow H, (x_i)_{1 \leq i \leq n} \mapsto \sum_{i=1}^n T_i x_i
\]
is a contraction. A dilation result of M\"uller and Vasilescu \cite{MV}, extended by Arveson \cite{Arv}, shows
that a tuple $T \in L(H)^n$ is a row contraction if and only if $T$ is up to unitary equivalence a compression of the
direct sum $M_z \oplus U \in L(H(\mathbb B,\mathcal D) \oplus K)^n$  of a Drury-Arveson shift $M_z \in L(H(\mathbb B,\mathcal D))^n$ and a
spherical unitary $U \in L(K)^n$ to one of its co-invariant subspaces. More precisely, let $\mathbb B \subset  \mathbb C^n$ be the open
Euclidean unit ball. Then by definition $M_z = (M_{z_1}, \ldots , M_{z_n}) \in L(H(\mathbb B,\mathcal D))^n$
is the tuple of multiplication operators with the coordinate functions on the $\mathcal D$-valued analytic functional Hilbert space 
$H(\mathbb B,\mathcal D)$, with a suitable Hilbert space $\mathcal D$, given by the reproducing kernel
\[
K: \mathbb B \times \mathbb B \rightarrow L(\mathcal D), (z,w) \mapsto \frac{1_{\mathcal D}}{1 -\langle z,w \rangle},
\]
while a spherical unitary is a commuting tuple $U = (U_1, \ldots, U_n) \in L(K)^n$ of normal operators with $\sum^n_{i=1} U_i U^*_i  = 1_K$.
The same condition, but without the spherical unitary part $U \in L(K)^n$, characterizes precisely the row contractions $T \in L(H)^n$ 
which are pure in the sense that they satisfy a $C_{\cdot 0}$-condition of the form
\[
{\rm SOT}-\lim_{k \rightarrow \infty} \sigma_T^k(1_H) = 0,
\]
where $\sigma_T: L(H) \rightarrow L(H)$ is the linear map defined by $\sigma_T(X) = \sum^n_{i=1} T_i X T^*_i$.

Let $D_T = (1_H - T^*T)^{1/2} \in L(H^n)$, $D_{T^*} = (1_H - TT^*)^{1/2} \in L(H)$ be the defect operators of $T$ and let
$\mathcal D_T = \overline{D_T H^n}$, $\mathcal D_{T^*} = \overline{D_{T^*} H}$ be the defect spaces of $T$, where
$T \in L(H^n,H)$ is regarded as a row operator and $T^* \in L(H,H^n)$ as a column operator. One way to prove the dilation
result mentioned above in the case of a pure row contraction $T \in L(H)^n$ is to show that the characteristic function
$\theta_T: \mathbb B \rightarrow L(\mathcal D_T,\mathcal D_{T^*})$,
\[
\theta_T(z) = -T +  D_{T^*}(1_H - ZT^*)ZD_T
\]
defines a partially isometric multiplier from $H(\mathbb B,\mathcal D_T)$ to $H(\mathbb B,\mathcal D_{T^*})$ and that 
$T$ is unitarily equivalent to the compression of $M_z \in L(H(\mathbb B,\mathcal D_{T^*}))^n$ to the co-invariant subspace
\[
\mathbb H_T = H(\mathbb B,\mathcal D_{T^*}) \ominus \theta_T H(\mathbb B,\mathcal D_T).
\]

Let $m > 0$ be a positive integer and let $H_m(\mathbb B,\mathcal D)$ be the $\mathcal D$-valued functional Hilbert space
given by the reproducing kernel
\[
K_m: \mathbb B \times \mathbb B \rightarrow L(\mathcal D), (z,w) \mapsto \frac{1_{\mathcal D}}{(1 -\langle z,w \rangle)^m}.
\]
Then the corresponding multiplication tuple $M_z \in L(H_m(\mathbb B,\mathcal D))^n$ plays the role of a model tuple for
a class of commuting Hilbert-space tuples $T \in L(H)^n$ satisfying suitable higher order positivity conditions. To be more precise,
the tuple $T$ is a row contraction if and only if $(1 - \sigma_T)(1_H) \geq 0$. A commuting tuple $T \in L(H)^n$ is called a 
row-$m$-hypercontraction, or simply an $m$-hypercontraction, if
\[
(1 - \sigma_T)^k(1_H) \geq 0 \quad (k = 1, \ldots ,m).
\]
In the cited paper \cite{MV} of M\"uller and Vasilescu it is shown that a commuting tuple $T \in L(H)^n$ is an
$m$-hypercontraction if and only if, up to unitary equivalence, it is the compression to a co-invariant subspace of a direct sum
$M_z \oplus U \in L(H_m(\mathbb B,\mathcal D) \oplus K)^n$ of the multiplication tuple $M_z = (M_{z_1}, \ldots , M_{z_n}) \in L(H_m(\mathbb B,\mathcal D))^n$
and a spherical unitary $U \in L(K)^n$. The situation in the general case $m \geq 1$ is completely analogous to the particular case of row conctractions
(that is, $m = 1$), except that in the proof given by M\"uller and Vasilecu no characteristic function is constructed. Indeed, up to now, no reasonable
definition of a characteristic function $\theta_T$ for $m$-hypercontractions seems to be known, except for some partial one-dimensional results
due to Ball and Bolotnikov \cite{BB1,BB2}. It is one of the purposes of the present paper
to suggest a possible definition of a characteristic function for $m$-hypercontractions.

In the classical Sz.-Nagy-Foias theory \cite{Sz-N}, which corresponds to the choices $n = 1$ and $m = 1$, the characteristic function $\theta_T$ of a pure
contraction $T \in L(H)$ induces an isometric multiplier $M_{\theta_T}: H^2(\mathbb D,\mathcal D_T) \rightarrow H^2(\mathbb D,\mathcal D_{T^*})$
between Hardy spaces.
In particular, the induced map $\theta_T: \mathcal D_T \rightarrow H^2(\mathbb D,\mathcal D_{T^*})$, $x \mapsto \theta_T x$, is isometric and
the closed subspace $\theta_T \mathcal D_T \subset H^2(\mathbb D,\mathcal D_{T^*})$ satisfies
\[
\theta_T \mathcal D_T \perp M^k_z(\theta_T \mathcal D_T) \quad \quad (k \geq 1).
\]
Closed subspaces with this property are usually called wandering subspaces for $M_z \in L(H^2(\mathbb D,\mathcal D_{T^*}))^n$.
In an effort to extend the notion of inner functions from the Hardy space setting to the case of Bergman spaces, Hedenmalm \cite{He} called a
function $f$ in the Bergman space $L^2_a(\mathbb D)$ on the unit disc Bergman inner if 
\[
\int_{\mathbb D} (|f|^2 - 1) z^k {\rm d}z = 0 \quad \quad (k \geq 0).
\]
If one regards $f: \mathbb D \rightarrow \mathbb C \cong L(\mathbb C)$ as an operator-valued function, then the preceding definition 
means precisely that the map $\mathbb C \rightarrow L^2_a(\mathbb D)$, $\alpha \mapsto f \alpha$, is isometric and 
\[
f\mathbb C \perp M^k_z(f \mathbb C) \quad \quad (k \geq 1).
\]
More generally, an operator-valued function $W: \mathbb B \rightarrow L(\mathcal E_*,\mathcal E)$ is called $K_m$-inner if the induced map
$\mathcal E_* \rightarrow H_m(\mathbb B,\mathcal E)$, $x \mapsto Wx$, is isometric and if its range satisfies the above orthogonality relations
\[
W \mathcal E_* \perp M_z^k(W \mathcal E_*) \quad \quad (k \geq 1).
\]
In the one-variable case, it was observed by Olofsson \cite{O1,O2} that each $K_m$-inner function
$W: \mathbb D \rightarrow L(\mathcal E_*,\mathcal E)$ on the unit disc admits a realization as a transfer function similar to the
one for characteristic functions and that, conversely, each function that admits such a realization defines a $K_m$-inner function.
Up to now no such realization formula is known for the higher-dimensional case. A second purpose of this note is to extend the
cited results of Olofsson from the case of the unit disc to the unit ball, and at the same time, to associate with each pure
$m$-hypercontraction $T \in L(H)^n$ a canonical $K_m$-inner function which is closely related to the characteristic
function $\theta_T$ of $T$.

In Section 2 we associate with each pure $m$-hypercontraction $T \in L(H)^n$ a canonical $K_m$-inner function
$W_T: \mathbb B \rightarrow L(\tilde{\mathcal D},\mathcal D)$ such that the space $W \tilde{\mathcal D} \subset H_m(\mathbb B,\mathcal D)$
is the wandering subspace of the invariant subspace $M = (j_T H)^{\perp} \in {\rm Lat}(M_z,H_m(\mathbb B,\mathcal D))$, where 
$j_T: H \rightarrow H_m(\mathbb B,\mathcal D)$ is the isometric intertwiner for $T^* \in L(H)^n$ and $M^*_z \in L(H_m(\mathbb B,\mathcal D))^n$
which yields the canonical functional model for $T$. In Section 3 we show that the $K_m$-inner functions 
$W: \mathbb B \rightarrow L(\mathcal E_*,\mathcal E)$ are precisely the operator-valued functions admitting a suitable transfer function
realization given by a pure $m$-hypercontraction $T \in L(H)^n$. In Section 4 we suggest a possible definition of a characteristic
function for $m$-hypercontractions.

\vspace{1cm }

\centerline{\textbf{\S2\, Wandering subspaces}}

Let $T=(T_1,\ldots,T_n)\in L(H)^n$ be a commuting tuple of bounded linear operators an a complex Hilbert space $H$. We use the bounded operator
\[
\sigma_T:L(H)\rightarrow L(H), X\mapsto \sum^n_{i=1}T_i X T^\ast_i
\]
to define the $k$-th order defect operators
\[
\Delta^{(k)}_T=(1-\sigma_T)^k(1_H)=\sum^k_{j=0}(-1)^j\binom{k}{j} \sum_{|\alpha|=j}\gamma_\alpha T^\alpha T^{\ast \alpha} \quad (k \in \mathbb N),
\]
where $\gamma_{\alpha} = |\alpha |!/\alpha !$ for $\alpha \in \mathbb N^n$.
We call $T$ a row-$m$-hypercontraction or simply an $m$-hypercontraction if the first and $m$-th order defect operators of $T$ are positive, that is,
\[
\Delta^{(1)}_T \geq 0 \mbox{ and }\Delta^{(m)}_T \geq 0.
\]
Commuting multioperators satisfying positivity conditions of this type have been studied for instance in \cite{MV}.\\
Let $T \in L(H)^n$ be an $m$-hypercontraction. It is well known (see Lemma 2 in \cite{MV}) that in this case even all defect operators up to order $m$ are positive, that is,
$\Delta^{(k)}_T\geq 0$ for $k=0,\ldots,m$. For $m=1$, an $m$-hypercontraction is usually called a row contraction. An $m$-hypercontraction $T \in L(H)^n$ is said to be pure or of
class $C_{\cdot 0}$ if
\[
{\rm SOT}-\lim_{k \rightarrow \infty}\sigma^k_T(1_H)=0.
\]
In \cite{MV} it was shown that each $m$-hypercontraction $T\in L(H)^n$ of class $C_{\cdot 0}$ is, up to unitary equivalence, 
the compression of a standard weighted shift to a co-invariant subspace.
To formulate this result in more detail, we need some additional notation. For positive integers $\ell, m \geq 1$ and each multiindex $\alpha \in \mathbb N^n$, we define
\[
\rho_\ell(\alpha)=\frac{(\ell + |\alpha|-1)!}{\alpha!(\ell - 1)!}.
\]
Let $T \in L(H)^m$ be an $m$-hypercontraction of class $C_{\cdot 0}$. We define $C = (\Delta^{(m)}_T)^{1/2}\in L(H)$ and $\mathcal D=\overline{CH}\subset H$. In
\cite{MV} (Theorem 9 and its proof) it was shown that the map
\[
H \rightarrow \ell^2(\mathbb N^n,\mathcal D), x \mapsto (\rho_m(\alpha)^{1/2}CT^{\ast\alpha}x)_{\alpha \in \mathbb N^n}
\]
is an isometry which intertwines the adjoint tuple
$T^\ast=(T^\ast_1,\ldots,T^\ast_n)\in L(H)^n$ and the backward shift tuple $S^{(m)}\in L(\ell^2(\mathbb N^n,\mathcal D))^n$ defined by
\[
(S^{(m)}_jx)(\alpha)=(\frac{\rho_m(\alpha)}{\rho_m(\alpha+e_j)})^{1/2}x(\alpha + e_j)\qquad (\alpha \in \mathbb N^n,j=1,\ldots,n)
\]
componentwise. Up to unitary equivalence the adjoint tuple $S^{(m)\ast}$ acts as the multiplication tuple $M_z=(M_{z_1},\ldots,M_{z_n})$ 
with the coordinate functions on a standard
analytic functional Hilbert space on the open unit ball. More precisely, given an integer $\ell \geq 0$ and a complex Hilbert space $\mathcal E$, we denote by $H_\ell(\mathbb B,\mathcal E)$ the $\mathcal E$-valued 
analytic functional Hilbert space with reproducing kernel
\[
K_\ell=K^\mathcal E_\ell:\mathbb B \times \mathbb B \rightarrow L(\mathcal E),K_\ell(z,w)=\frac{1_\mathcal E}{(1-\langle z,w\rangle)^\ell}
\]
on the Euclidean unit ball $\mathbb B=\lbrace z \in \mathbb C^n;\| z \| < 1\rbrace$. Then $H_0(\mathbb B,\mathcal E)=\mathcal E$ and
\[
H_\ell(\mathbb B,\mathcal E) = \{ f = \sum_{\alpha \in \mathbb N^n} f_\alpha z^\alpha \in \mathcal O(\mathbb B,\mathcal E); \| f \|^2=\sum_{\alpha \in \mathbb N^n}\frac{\| 
f_\alpha\|^2}{\rho_\ell(\alpha)}< \infty \}
\]
for $\ell \geq 1$. The spaces $H_1(\mathbb B,\mathcal E),H_n(\mathbb B,\mathcal E)$ and $H_{n+1}(\mathbb B,\mathcal E)$ are 
the $\mathcal E$-valued Drury-Arveson space $H^2_n(\mathbb B,\mathcal E)$, the
Hardy space $H_2(\mathbb B,\mathcal E)$ and the (unweighted) Bergman space $L^2_a(\mathbb B,\mathcal E)$, respectively.\\
Modulo the unitary operator
\[
U:\ell^2(\mathbb N^n,\mathcal D)\rightarrow H_m(\mathbb B,\mathcal D),(x_\alpha)\mapsto \sum_{\alpha \in \mathbb N^n}\rho_m(\alpha)^{1/2}x_\alpha z^\alpha
\]
the forward shift tuple $S^{(m)\ast}\in L(\ell^2(\mathbb N^n,\mathcal D))^n$ is unitarily equivalent to the multiplication tuple $M_z\in L(H_m(\mathbb B,\mathcal D))^n$. Indeed (see Equation (13) in \cite{MV})
\[
US^{(m)\ast}_j(x_\alpha)=\sum_{\substack{\alpha \in \mathbb N^n\\ \alpha_j \geq 1}} \rho_m(\alpha-e_j)^{1/2}x_{\alpha-e_j}z^\alpha =M_{z_j}U(x_\alpha).
\]
Hence the map $j:H\rightarrow H_m(\mathbb B,\mathcal D)$,
\[
j(x)=\sum_{\alpha \in \mathbb N^n}\rho_m(\alpha)(CT^{\ast\alpha}x)z^\alpha=C(1_H-ZT^\ast)^{-m}x
\]
defines an isometry intertwining the tuples $T^\ast \in L(H)^n$ and $M^\ast_z\in L(H_m(\mathbb B,\mathcal D))^n$ componentwise. 
Here the right-hand side is regarded as a function in $z \in \mathbb B$, we write $Z: H^n \rightarrow H, (h_i) \mapsto \sum^n_{i=1} z_i h_i$,
$T: H^n \rightarrow H, (h_i) \mapsto \sum^n_{i=1} T_i h_i$, for the row operators and $T^*: H \rightarrow H^n, h \mapsto (T^*_i h)^n_{i=1}$, for
the associated column operator. Following \cite{K} we call any such interzwining map an $m$-dilation for $T$.
It follows that the closed subspace
\[
M=H_m(\mathbb B,\mathcal D)\ominus ({\rm Im}\ j)\subset H_m(\mathbb B,\mathcal D)
\]
is invariant for the multiplication tuple $M_z \in L(H_m(\mathbb B,\mathcal D))^n$. Our aim in this section is to derive an explicit description of its wandering subspace
\[
W(M) = M\ominus (\sum^n_{i=1}M_{z_i}M)
\]
which extends corresponding one-variable results from \cite{O1}.\\
Our main tool will be the matrix operator $M^\ast_zM_z=(M^\ast_{z_i}M_{z_j})_{1\leq i,j\leq n} \in L(H_m(\mathbb B,\mathcal D)^n)$. Since the row operator $M_z:H_m(\mathbb B,\mathcal D)^n\rightarrow H_m(\mathbb B,\mathcal D)$ has closed range
\[
M_zH_m(\mathbb B,\mathcal D)^n=\lbrace f \in H_m(\mathbb B,\mathcal D);f(0)=0\rbrace,
\]
the operator $M^\ast_zM_z:{\rm Im}M^\ast_z\rightarrow {\rm Im}M^\ast_z$ is invertible. We use the notation $(M^\ast_zM_z)^{-1}$ for its inverse. As a first step we show that the action of this inverse
can be described using the operator $\delta : H_m(\mathbb B,\mathcal D)\rightarrow H_m(\mathbb B,\mathcal D)$,
\[
\delta (\sum^\infty_{k=0} \sum_{| \alpha|=k}f_\alpha z^\alpha)  = f_0 + \sum^\infty_{k=1} (\frac{m+k-1}{k}\sum_{| \alpha |=k}f_\alpha z^\alpha).
\]
Note that $\delta$ is a diagonal operator with respect to the orthogonal decomposition
\[
H_m(\mathbb B,\mathcal D) = \oplus^{\infty}_{k=0} \mathbb H_k(\mathcal D)
\]
of $H_m(\mathbb B,\mathcal D)$ into the subspaces $\mathbb H_k(\mathcal D)$ consisting of all homogeneous polynomials 
of degree $k$ with coefficients in $\mathcal D$.

%Lemma 2.1
\begin{lem}
\label{inverse}
For $f \in H_m(\mathbb B,\mathcal D)$,
\[
(M^\ast_z M_z)^{-1}(M^\ast_z f)=M^\ast_z \delta f.
\]
Im particular, the row operator
\[
H_m(\mathbb B,\mathcal D)^n\stackrel{\delta M_{z}}{\longrightarrow}H_m(\mathbb B,\mathcal D)
\]
defines a continuous linear extension of the operator
\[
M_z(M^\ast_zM_z)^{-1}: {\rm Im}M^\ast_z \rightarrow H_m(\mathbb B,\mathcal D).
\]
\end{lem}

\proof
An elementary calculation shows that
\[
M^\ast_{z_i}g=\sum_{\alpha \in \mathbb N^n}\frac{\rho_m(\alpha)}{\rho_m(\alpha + e_i)}g_{\alpha + e_i}z^\alpha=\sum_{\alpha \in \mathbb N^n}\frac{\alpha_i+1}{m+|\alpha|}g_{\alpha + e_i}z^\alpha
\]
for $g=\sum_{\alpha \in \mathbb N^n}g_\alpha z^\alpha \in H_m(\mathbb B,\mathcal D)$ and $i=1,\ldots,n$. Let $f \in H_m(\mathbb B,\mathcal D)$ 
be arbitrary. To prove the first assertion we may suppose
that $f(0)=0$. It is easy to see and well known (see e.g. Section 3 in \cite{AK}) that in this case $f=\sum^n_{i=1}z_if_i$ with
\[
f_i = \sum_{\alpha \in \mathbb N^n}\frac{\alpha_i+1}{|\alpha|+1}f_{\alpha + e_i}z^\alpha = 
\sum_{\alpha \in \mathbb N^n}\frac{\alpha_i+1}{m+|\alpha|}(\frac{m+|\alpha|}{| \alpha|+1}f_
{\alpha + e_i})z^\alpha=M^\ast_{z_i} \delta f
\]
for $i=1,\ldots,n$. But then
\[
(M^\ast_zM_z)^{-1}M^\ast_z f = (M^\ast_zM_z)^{-1}(M^\ast_zM_z)(M^\ast_z \delta f) = M^\ast_z \delta f.
\]
Since $M_z M^*_z \in L(H_m(\mathbb B,\mathcal D))$ is a diagonal operator, we have  
$M_z(M^\ast_zM_z)^{-1}M^\ast_z = M_zM^\ast_z \delta = \delta M_zM^\ast_z$. Hence also the second assertion follows.
\proofend

Let $P_\mathcal D \in L(H_m(\mathbb B,\mathcal D))$ be the orthogonal projection onto the closed subspace consisting of all constant functions. Then ${\rm Ker}\ P_\mathcal D=\lbrace f \in H_m(\mathbb B,\mathcal D);f(0)=0\rbrace$. 

%Lemma 2.2
\begin{lem}
\label{projection}
The orthogonal projection $P_\mathcal D$ acts as
\[
P_\mathcal D=1_{H_m(\mathbb B,\mathcal D)}-\sum^{m-1}_{j=0}(-1)^j\binom{m}{j+1}\sum_{| \alpha|=j+1}\gamma_\alpha M^\alpha_zM^{\ast \alpha}_z.
\]
\end{lem}

\proof
Let us write
\[
L = (L_{M_{z_1}},\ldots,L_{M_{z_n}}),\ R = (R_{M^\ast_{z_1}},\ldots, R_{M^\ast_{z_n}}) \in L(L(H_m(\mathbb B,\mathcal D)))^n
\]
for the tuples consisting of the left and right multiplication operators
\[
L_{M_{z_i}}(X)=M_{z_i}X,\ R_{M^\ast_{z_i}}(X)=XM^\ast_{z_i}\qquad (X\in L(H_m(\mathbb B,\mathcal D))).
\]
By Lemma 1.2 in \cite{AE} the projection $P_\mathcal D$ is given by
\[
P_\mathcal D=\frac{1}{C}(L,R)(1_{H_m(\mathbb B,\mathcal D)}),
\]
where $C(z,w)=K_m(z,\overline{w})$ and $(1/C)(z,w)=(1-\langle z,\overline{w}\rangle)^m$ is regarded as the analytic polynomial
\[
(1-\langle z,\overline{w}\rangle)^m=\sum^m_{j=0}(-1)^j\binom{m}{j}\sum_{|\alpha|=j}\gamma_\alpha z^\alpha w^\alpha.
\]
Thus the assertion follows.
\proofend

A more elementary, alternative representation of the orthogonal projection $P_\mathcal D$ follows from the observation that 
the operator $M_z(M^\ast_zM_z)^{-1}M^\ast_z \in L(H_m(\mathbb B,\mathcal D))$ 
is the orthogonal projection onto the space ${\rm Im}\ M_z=\lbrace f \in H_m(\mathbb B,\mathcal D);\ f(0)=0\rbrace ={\rm Ker}\ P_\mathcal D$. 
The above preparations allow us to deduce a first description of the wandering subspace $W(M) = M \ominus (\sum^n_{i=1} M_{z_i} M)$.

%Theorem 2.3
\begin{thm}
\label{wandering1}
A function $f \in H_m(\mathbb B,\mathcal D)$ belongs to the wandering subspace $W(M)$ of 
$M = ({\rm Im}j)^\perp \in {\rm Lat}(M_z,H_m(\mathbb B,\mathcal D))$ if and only if
\[
f = f_0 + M_z(M^\ast_zM_z)^{-1}(jx_i)^n_{i=1}
\]
for some vectors $f_0 \in \mathcal D,x_1,\ldots,x_m \in H$ with $(jx_i)^n_{i=1} \in M^\ast_zH_m(\mathbb B,\mathcal D)$ and
\[
Cf_0 + T \big( \sum^{m-1}_{j=0}(-1)^j\binom{m}{j+1}\sum_{|\alpha|=j}\gamma_\alpha T^\alpha T^{\ast\alpha}x_i)\big)^n_{i=1}=0.
\]
In this case, the identity $(jx_i)^n_{i=1} = M^*_z f$ holds.
\end{thm}

\proof
Since $M = \ker j^* \in {\rm Lat}(M_z)$, it follows that a function
$f\in H_m(\mathbb B,\mathcal D)$ belongs to the space $W(M) = M\ominus(\sum^n_{i=1}z_iM)$ if and only if 
$j^\ast f = 0$ and $(1_{H_m(\mathbb B,\mathcal D)}-jj^\ast)M^\ast_{z_i}f = 0$ for $i = 1,\ldots,n$. 
Suppose first that $f$ satisfies these conditions. Then $(x_i)^n_{i=1} = (j^{\ast}M^\ast_{z_i}f)^n_{i=1}$ 
defines a tuple in $H^n$ with $(jx_i)^n_{i=1} = M^\ast_zf$ such that
\[
f = f(0)+M_z(M^\ast_zM_z)^{-1}M^\ast_zf = f(0) + M_z(M^\ast_zM_z)^{-1}(jx_i)^n_{i=1}.
\]
Using the definition of the isometry $j: H\rightarrow H_m(\mathbb B,\mathcal D)$ we find that
\[
\langle y,j^\ast x\rangle=\langle jy,x\rangle=\langle j(y)(0),x\rangle=\langle Cy,x\rangle=\langle y,Cx\rangle
\]
for all $y \in H$ and $x \in \mathcal D$ and hence that $j^\ast x = Cx$ for $x \in \mathcal D$ regarded as a constant 
function in $H_m(\mathbb B,\mathcal D)$. Using Lemma \ref{inverse} and the intertwining properties of $j$ we obtain that
\[
0 = j^* f = Cf(0) + T(\oplus j^*)(M^*_zM_z)^{-1}M^*_z f = Cf(0) + T(\oplus j^*)M^*_z \delta f.
\]
A straightforward calulation, or the results in Section 3 of \cite{AK}, show that
\[
\langle z,w\rangle (M^*_{z_i} \delta K_m(\cdot,w)x)(z)=\overline{w}_i(K_m(z,w)x-x)
\]
for $z,w \in \mathbb B$, $x \in \mathcal D$ and $i = 1, \ldots ,n$. Using Lemma \ref{projection} we find that
\[
\langle z,w \rangle \Bigl( \sum^{m-1}_{j=0}(-1)^j\binom{m}{j+1}\sum_{|\alpha |=j}\gamma_\alpha M^\alpha_z M^{\ast\alpha}_z M^\ast_{z_i}K_m(\cdot,w)x \Bigr) (z)
\]
\[
= \overline{w}_i\sum^{m-1}_{j=0}(-1)^j\binom{m}{j+1}\langle z,w\rangle^{j+1}K_m(z,w)x
\]
\[
= \overline{w}_i \Big( \sum^{m-1}_{j=0}(-1)^j\binom{m}{j+1}\sum_{|\alpha|=j+1}\gamma_\alpha M^\alpha_zM^{\ast\alpha}_zK_m(\cdot,w)x \Bigr)(z)
\]
\[
= \overline{w}_i((1_{H_m(\mathbb B,\mathcal D)}-P_\mathcal D)K_m(\cdot,w)x)(z)
=\overline{w}_i(K_m(z,w)x-x)
\]
for $z,w \in \mathbb B$, $x\in \mathcal D$ and $i = 1 \ldots ,n$. By comparing the previous two results and using the fact
that the closed linear span of the functions $K_m(\cdot,w)x$ $(w \in \mathbb B, x \in \mathcal D)$ is of of
$H_m(\mathbb B, \mathcal D)$, we obtain
\begin{align*}
0& = Cf(0) + T(\oplus j^*) \Big( \sum^{m-1}_{j=0}(-1)^j\binom{m}{j+1}\sum_{|\alpha|=j}\gamma_\alpha M^\alpha_z M^{* \alpha}_z M^*_{z_i}f \Bigr)^n_{i=1}\\
&=Cf(0)+T(\sum^{m-1}_{j=0}(-1)^j\binom{m}{j+1}\sum_{|\alpha|=j}\gamma_\alpha T^\alpha T^{\ast\alpha}x_i)^n_{i=1}.
\end{align*}
Conversely, suppose that $f \in H_m(\mathbb B,\mathcal D)$ is a function that has a representation as described in 
Theorem \ref{wandering1}. We show that $f\in W(M)$. Note first that
\[
M^\ast_zf=M^\ast_z(f_0+M_z(M^\ast_zM_z)^{-1}(jx_i)^n_{i=1})=(jx_i)^n_{i=1}
\]
and therefore $jj^\ast M^\ast_{z_i}f=jx_i=M^\ast_{z_i} f$ for $i=1,\ldots,n$. Then
exactly as in the first part of the proof it follows that
\begin{align*}
j^\ast f&=Cf_0+j^\ast M_z(M^\ast_zM_z)^{-1}M^\ast_zf\\
&=Cf_0+T(\sum^{m-1}_{j=0}(-1)^j\binom{m}{j+1}\sum_{|\alpha|=j}\gamma_\alpha T^\alpha T^{\ast\alpha}x_i)^n_{i=1}.
\end{align*}
Thus also the reverse implication follows.
\proofend

For future use, note that in the above proof we deduced the formula
\[
M^*_{z_i} \delta = \Big( \sum^{m-1}_{j=0}(-1)^j\binom{m}{j+1}\sum_{|\alpha|=j}\gamma_\alpha M^\alpha_z M_z^{* \alpha} \Bigr) M^*_{z_i} \quad (1 \leq i \leq n).
\] 

The defect operators $\Delta^{(k)}_T=(1-\sigma_T)^k(1_H)$ can be used to rewrite the conditions used in Theorem \ref{wandering1} to characterize the functions in $W(M)$.

%Lemma 2.4
\begin{lem}
\label{defect}
Let $T\in L(H)^n$ be an $m$-hypercontraction. Then
\[
\sum^{m-1}_{k=0}\Delta^{(k)}_T=\sum^{m-1}_{j=0}(-1)^j\binom{m}{j+1}\sum_{|\alpha|=j}\gamma_\alpha T^\alpha T^{\ast\alpha}.
\]
\end{lem}

\proof
Define
\[
\Sigma_k = \sum^{k-1}_{j=0}(-1)^j \binom{k}{j+1}\sum_{|\alpha|=j}\gamma_\alpha T^\alpha T^{\ast\alpha}\qquad (k=1,\ldots,m)
\]
and $\Sigma_0=0$. Then
\[
\Delta^{(k)}_T=\sum^k_{j=0}(-1)^j(\binom{k+1}{j+1}-\binom{k}{j+1})\sum_{|\alpha|=j}\gamma_\alpha T^\alpha T^{\ast\alpha}=\Sigma_{k+1}-\Sigma_k
\]
for $k=0,\ldots,m-1$. Hence $\sum^{m-1}_{k=0}\Delta^{(k)}_T=\Sigma_m$.
\proofend

An inspection of the proof of Theorem \ref{wandering1} allows us to calculate the norms of the functions $f \in W(M)$ in terms of the data $f_0,x_1,\ldots,x_m$ occurring in their standard representation.

%Lemma 2.5
\begin{lem}
\label{norm}
Let $T \in L(H)^n$ be an $m$-hypercontraction of class $C_{\cdot 0}$ and let
\[
f=f_0+M_z(M^\ast_zM_z)^{-1}(jx_i)^n_{i=1}
\]
be a representation of a function $f \in W(M)$ as in Theorem \ref{wandering1}. Then
\[
\| f \|^2=\| f_0\|^2+\sum^{m-1}_{j=0}(-1)^j\binom{m}{j+1}\sum_{|\alpha|=j}\gamma_\alpha(\sum^n_{i=1}\| T^{\ast\alpha}x_i\|^2).
\]
\end{lem}

\proof
Since ${\rm Im}M_z\subset H_m(\mathbb B,\mathcal D)\ominus \mathcal D$, we obtain
\[
\| f \|^2=\| f_0\|^2+\| M_z(M^\ast_zM_z)^{-1}(jx_i)^n_{i=1}\|^2.
\]
As observed in Theorem \ref{wandering1} $(jx_i)^n_{i=1} = M^\ast_z f$. Exactly as in the proof of Theorem \ref{wandering1}, 
it follows that the second term in the above sum is given by

\[
\langle(\oplus j^\ast)(M^\ast_zM_z)^{-1}M^\ast_z f,(x_i)^n_{i=1}\rangle
= \langle(\oplus j^\ast)M^\ast_z \delta f,(x_i)^n_{i=1}\rangle
\]
\[
= \langle(\sum^{m-1}_{j=0}(-1)^j\binom{m}{j+1}\sum_{|\alpha |=j}\gamma_\alpha T^\alpha T^{\ast\alpha}x_i)^n_{i=1},(x_i)^n_{i=1}\rangle.
\]

Thus the assertion follows.
\proofend

Let $T \in L(H)^n$ be an $m$-hypercontraction of class $C_{\cdot 0}$. Define $D_{T,k}=(\Delta^{(k)}_T)^{1/2}\in L(H)$ for $k=0,\ldots,m-1$. Then
\[
(x,y)=\sum^{m-1}_{k=0}\langle \Delta^{(k)}_Tx,y\rangle=\sum^{m-1}_{k=0}\langle D_{T,k}x,D_{T,k}y\rangle
\]
defines a scalar product on $H$ such that the induced norm $\| x \|_m=(x,x)^{1/2}$ is equivalent to the original norm of $H$. We write $\tilde{H}$ for $H$ equipped with this norm. Then $I_m: H\rightarrow \tilde{H}, x \mapsto x$, defines an invertible bounded operator such that
\[
\langle I^\ast_mx,y\rangle = \langle\sum^{m-1}_{k=0}  \Delta^{(k)}_T x,y\rangle\qquad (x \in \tilde{H}, y \in H).
\]
Hence $I^\ast_m: \tilde{H}\rightarrow H$ acts as $I^\ast_m x=\sum^{m-1}_{k=0}\Delta^{(k)}_Tx$. We write $\tilde{T}=(\tilde{T_1},\ldots,\tilde{T_n}):\tilde{H}^n\rightarrow H$ for the row operator with 
components $\tilde{T_i}=T_i\circ I^\ast_m \in L(\tilde{H},H)$. Then
\[
\tilde{T}\tilde{T}^\ast=\sum^n_{i=1}T_i(I^\ast_m I_m)T^\ast_i=\sigma_T(\sum^{m-1}_{k=0}(1-\sigma_T)^k(1_H))=1_H-\Delta^{(m)}_T
\]
and hence $\tilde{T}$ is a contraction. As in \cite{O1} we shall use its defect operators
\begin{align*}
D_{\tilde{T}^\ast}&=(1_H-\tilde{T}\tilde{T}^\ast)^{1/2}=(\Delta^{(m)}_T)^{1/2}=C\in L(H),\\
D_{\tilde{T}}&=(1_{\tilde{H}^n}-\tilde{T}^\ast\tilde{T})^{1/2}\in L(\tilde{H}^n)
\end{align*}
and defect spaces $\mathcal D_{\tilde{T}^\ast}=\overline{D_{\tilde{T}^\ast} H}=\mathcal D,\mathcal D_{\tilde{T}}=\overline{D_{\tilde{T}}(\tilde{H}^n)}$ to deduce an alternative description of the wandering subspace
$W(M)$. Exactly as in the single-variable theory of contractions it follows that $\tilde{T} D_{\tilde{T}}=D_{\tilde{T}^\ast}\tilde{T}$ and that
\[
U=\biggl(
\begin{array}{c|c}
\tilde{T}&D_{\tilde{T}^\ast}\\
\hline
D_{\tilde{T}}&-\tilde{T}^\ast
\end{array}
\biggr)
:\tilde{H}^n\oplus \mathcal D_{\tilde{T}^*} \rightarrow H\oplus \mathcal D_{\tilde{T}}
\]
is a well-defined unitary operator. In the following we construct an analytically parametrized family of operators $W_T(z)\in L(\tilde{\mathcal D},\mathcal D)$ on the closed subspace
\[
\tilde{\mathcal D}=\lbrace y \in \mathcal D_{\tilde{T}}; (\oplus j I^{-1}_m)D_{\tilde{T}}y\in M^\ast_zH_m(\mathbb B,\mathcal D)\rbrace \subset \mathcal D_{\tilde{T}}
\]
such that
\[
W(M)=\lbrace W_T x;\ x \in \tilde{\mathcal D}\rbrace,
\]
where $W_Tx: \mathbb B \rightarrow \mathcal D$ acts as $(W_T x)(z) = W_T(z)x.$

%Lemma 2.6
\begin{lem} \label{parameter}
Let $T \in L(H)^n$ be an $m$-hypercontraction of class $C_{\cdot 0}$. Then a function $f \in H_m(\mathbb B,\mathcal D)$ belongs to the wandering subspace $W(M)$ of $M=({\rm Im}j)^\perp \in {\rm Lat}(M_z,H_m(\mathbb B,\mathcal D))$
if and only if there is a vector $y \in \tilde{\mathcal D}$ with
\[
f=-\tilde{T} y+M_z(M^\ast_zM_z)^{-1}(\oplus j I^{-1}_m)D_{\tilde{T}}y.
\]
In this case, we have $\| f \|^2 = \| y \|^2_{\tilde{H}^n}$.
\end{lem}

\proof
By Theorem \ref{wandering1} and Lemma \ref{defect} the space $W(M)$ consists precisely of all functions 
$f \in H_m(\mathbb B,\mathcal D)$ of the form
\[
f = f_0 + M_z(M^\ast_zM_z)^{-1}(jx_i)^n_{i=1}
\]
where $f_0\in \mathcal D,x_1,\ldots,x_n\in H$ are vectors with $(jx_i)^n_{i=1}\in M^\ast_zH_m(\mathbb B,\mathcal D)$ and
\[
\tilde{T}(I_mx_i)^n_{i=1}+D_{\tilde{T}^\ast}f_0=0.
\]
Then $y=D_{\tilde{T}}(I_mx_i)^n_{i=1}-\tilde{T}^\ast f_0\in \mathcal D_{\tilde{T}}$ is a vector with
\[
U\binom{(I_mx_i)}{f_0}=\binom{0}{y},
\]
or equivalently, with
\[
\binom{(I_mx_i)}{f_0} = U^\ast\binom{0}{y} = \binom{D_{\tilde{T}}y}{-\tilde{T}y}.
\]
But then $y \in \tilde{\mathcal D}$ and
\[
f = -\tilde{T}y + M_z(M^\ast_zM_z)^{-1}(\oplus j I^{-1}_m)D_{\tilde{T}}y.
\]
Conversely, if $f$ is a function of this form for some vector $y \in \tilde{\mathcal D}$, then by reversing the above arguments, one finds that
\[
f_0 = -\tilde{T}y\in \mathcal D,(x_i)^n_{i=1} = (\oplus I^{-1}_m)D_{\tilde{T}}y\in H^n
\]
are vectors with $(jx_i)^n_{i=1}\in M^\ast_zH_m(\mathbb B,\mathcal D)$ and
\[
f = f_0+M_z(M^\ast_zM_z)^{-1}(jx_i)^n_{i=1}.
\]
Hence $f\in W(M)$ by Theorem \ref{wandering1}. By Lemma \ref{defect} and Lemma \ref{norm},
\[
\| f \|^2 = \| f_0\|^2+\sum^n_{i=1}\| x_i\|^2_m = \| \tilde{T}y\|^2+\| D_{\tilde{T}}y\|^2_{\tilde{H}^n} = \| y\|^2_{\tilde{H}^n}.
\]
This observation completes the proof.
\proofend

We define an analytic operator-valued function $W_T:\mathbb B\rightarrow L(\tilde{\mathcal D},\mathcal D)$ by
\[
W_T(z)(x_i)^n_{i=1} = -T(\sum^{m-1}_{k=0}\Delta^{(k)}_TI^{-1}_mx_i)+C\sum^m_{k=1}(1_H-ZT^\ast)^{-k}I^{-1}_mZD_{\tilde{T}}(x_i)^n_{i=1}.
\]
Note that the first term in the defining sum on the right equals
\[
W_T(0)(x_i)^n_{i=1} = - \tilde{T}(x_i)^n_{i=1} \in \mathcal D.
\]

%Theorem 2.7
\begin{thm}
\label{wandering2}
Let $T \in L(H)^n$ be an $m$-hypercontraction of class $C_{\cdot 0}$. Then $W(M)=\lbrace W_Tx;x\in \tilde{\mathcal D}\rbrace$ and
\[
\| W_T x \|=\| x \|\qquad (x \in \tilde{\mathcal D}).
\]
\end{thm}

\proof
Exactly as in \cite{Ol} (proof of Lemma 1.1 and page 534) we use the formula
\[
\sum^m_{k=1}(1-z)^{-k}=\sum^\infty_{j=0}\rho_m(j+1)z^j\qquad (z\in \mathbb D)
\]
to obtain the representations
\[
\sum^m_{k=1}(1_H-ZT^\ast)^{-k}=\sum^\infty_{j=0}\rho_m(j+1)(ZT^\ast)^j=\sum^\infty_{j=0}\rho_m(j+1)\sum_{|\alpha|=j}(\gamma_\alpha T^{\ast\alpha})z^\alpha
\]
for $z \in \mathbb B$. A simple calculation shows that
\[
\rho_m(j+1)\ \gamma_\alpha=\frac{m+|\alpha|}{|\alpha|+1}\rho_m(\alpha)
\]
for $j \in \mathbb N$ and $\alpha \in \mathbb N^n$ with $|\alpha|=j$. Let $D \in L(H_m(\mathbb B,\mathcal D))$ be the diagonal operator used in Lemma \ref{inverse}. Then
\[
C \sum^m_{k=1}(1_H-ZT^\ast)^{-k}Z(x_i)^n_{i=1}
= \sum^n_{i=1}\sum^\infty_{j=0}\rho_m(j+1)(\sum_{|\alpha|=j}\gamma_\alpha CT^{\ast \alpha} x_iz^{\alpha + e_i})
\]
\[
= \sum^n_{i=1} \delta M_{z_i}(\sum_{\alpha \in \mathbb N^n}\rho_m(\alpha) CT^{\ast\alpha}x_i z^\alpha)
= (\sum^n_{i=1} \delta M_{z_i}jx_i)(z)
\]
for $z \in \mathbb B$ and $(x_i)^n_{i=1}\in H^n$. If $(x_i)^n_{i=1}\in H^n$ is a tuple with $(jx_i)^n_{i=1}\in {\rm Im}M^\ast_z$, then by Lemma \ref{inverse}
\[
M_z(M^\ast_zM_z)^{-1}(jx_i)^n_{i=1} = \sum^n_{i=1}DM_{z_i}(jx_i).
\]
Comparing the previous two results we obtain the identity
\[
C\sum^m_{k=1}(1_H-ZT^\ast)^{-k}I^{-1}_mZD_{\tilde{T}}x = M_z(M^\ast_zM_z)^{-1}(\oplus j I^{-1}_m)D_{\tilde{T}}x
\]
for every vector $x \in \tilde{\mathcal D}$. Thus an application of Lemma \ref{parameter} completes the proof.
\proofend

We briefly indicate an alternative description of the defect space $\tilde{\mathcal D}$. For a commuting tuple $S\in L(H)^n$, let 
$H^p(S,H) = {\rm Ker}\ \delta^p_S/{\rm Im}\ \delta^{p-1}_S$ $(p=0,\ldots,n)$ be the cohomology groups of its Koszul complex
(cf. Chapter 2 in \cite{EP})
\[
K^{\cdot}(S,H):0\rightarrow \Lambda^0(\sigma,H)\stackrel{\delta^0_S}{\longrightarrow} \Lambda^1(\sigma,H)\stackrel{\delta^{1}_S}{\longrightarrow}\ldots \stackrel{\delta^{n-1}_{S}}{\longrightarrow} 
\Lambda^n(\sigma,H)\rightarrow 0.
\]
 Then $H^p(M_z,H_m(\mathbb B,\mathcal D)) = 0$ for $p < n$ and $H^n(M_z,H_m(\mathbb B,\mathcal D))\cong \mathcal D$. 
By duality we obtain that $H^0(M^\ast_z,H_m(\mathbb B,\mathcal D))\cong \mathcal D$ and $H^p(M^\ast_z,H_m(\mathbb B,\mathcal D))=0$ for $p> 0$.

\begin{lem} \label{Koszul}
Let $T\in L(H)^n$ be an $m$-hypercontraction of class $C_{\cdot 0}$. Then
\[
\tilde{\mathcal D}=\lbrace y \in \mathcal D_{\tilde{T}};(\oplus I^{-1}_m)D_{\tilde{T}}y \in {\rm Ker}\ \delta^1_{T^\ast}\rbrace.
\]
\end{lem}

\proof
Let $y \in \mathcal D_{\tilde{T}}$ be arbitrary. We write $x=(x_i)^n_{i=1}=(\oplus I^{-1}_m)D_{\tilde{T}}y$ for $D_{\tilde{T}}y$ regarded as an element in $H^n$. If $y \in \tilde{\mathcal D}$, then
$(\oplus j)x=M^\ast_z f$ for some function $f \in H_m(\mathbb B,\mathcal D)$. Since
\[
j T^\ast_kx_i = M^\ast_{z_k}jx_i = M^\ast_{z_k}M^\ast_{z_i}f = j(T^\ast_ix_k)\qquad (i,k=1,\ldots,n),
\]
it follows that $x \in {\rm Ker}\ \delta^1_{T^\ast}$. Conversely, if $x \in {\rm Ker}\ \delta^1_{T^\ast}$, then 
\[
M^\ast_{z_k}jx_i=jT^\ast_kx_i=jT^\ast_ix_k=M^\ast_{z_i}jx_k\qquad (i,k=1,\ldots,n).
\]
Since $H^1(M^\ast_z,H_m(\mathbb B,\mathcal D)) = 0$, it follows that $jx \in M^*_z H_m(\mathbb B, \mathcal D)$ and hence that $y \in \tilde{\mathcal D}$.
\proofend

Since in the setting of Theorem \ref{wandering2} the space $W_T(\tilde{\mathcal D})=W(M)$ is a wandering subspace for $M_z|_M\in L(M)^n$, it follows that $W_T:\mathbb B\rightarrow L(\tilde{\mathcal D},\mathcal D)$
is an operator-valued analytic function with $W_T(x)\in H_m(\mathbb B,\mathcal D), \| W_Tx\|=\| x\|$ for all $x \in \tilde{\mathcal D}$ and
\[
W_T(\tilde{\mathcal D})\,\perp\, M^\alpha_z(W_T(\tilde{\mathcal D}))\mbox{ for all } \alpha \in \mathbb N^n\setminus \lbrace 0 \rbrace.
\]
This means precisely that $W_T:\mathbb B \rightarrow L(\tilde{\mathcal D},\mathcal D)$ is a $K_m$-inner function in the sense of \cite{BEKS}. As an application of Theorem 6.2 in \cite{BEKS} we obtain that $W_T$ is a contractive
multiplier from $H_1(\mathbb B,\tilde{\mathcal D})$ to $H_m(\mathbb B,\mathcal D)$.

\begin{cor}
\label{multiplier}
The operator-valued function $W_T:\mathbb B \rightarrow L(\tilde{\mathcal D},\mathcal D)$ induces a contractive multiplication operator
\[
H_1(\mathbb B,\tilde{\mathcal D})\rightarrow H_m(\mathbb B,\mathcal D),f\mapsto W_T f.
\]
\end{cor}

For $m = 1$, the set of $m$-hypercontractions $T\in L(H)^n$ of class $C_{\cdot 0}$ coincides with the class 
of all pure row contractions. For a pure row contraction $T\in L(H)^n$, the operator-valued function
\[
\theta_T:\mathbb B \rightarrow L(\mathcal D_T,\mathcal D_{T^\ast}), \theta_T(z)=-T+D_{T^\ast}(1_H-ZT^\ast)^{-1}ZD_T
\]
is called the {\rm characteristic function} of $T$. It is known \cite{BES1, BES2} that $\theta_T$ induces a partially isometric multiplication operator
\[
M_{\theta_T}:H_1(\mathbb B,\mathcal D_T)\rightarrow H_1(\mathbb B,\mathcal D_{T^\ast}), f \mapsto \theta_T f
\]
such that $M_{\theta_T}M^\ast_{\theta_T}+jj^\ast = 1_{H_1(\mathbb B,\mathcal D_{T^\ast})}$. In particular, $M=({\rm Im}j)^\perp={\rm Im}M_{\theta_T}$.

\begin{cor} \label{row}
Let $T\in L(H)^n$ be a pure row contraction. Then $W_T(z)=\theta_T(z)|_{\tilde{\mathcal D}}$ for all $z \in \mathbb B$ and the 
characteristic function $\theta_T$ induces a unitary operator
\[
\tilde{\mathcal D}\rightarrow W(M), x\mapsto M_{\theta_T}x.
\]
\end{cor}

\proof
In the particular case $m = 1$, our previously constructed spaces and operators reduce to 
$\tilde{H} = H,\tilde{T} = T$, $C = D_{T^*}$ and $\mathcal D_{\tilde{T}}=\mathcal D_T,\mathcal D=\mathcal D_{\tilde{T}^\ast}=\mathcal D_{T^\ast}$.
The domain of the operators $W_T(z)\ (z\in \mathbb B)$ is given by
\[
\tilde{\mathcal D}=\lbrace y \in \mathcal D_T;(\oplus j)D_Ty\in M^\ast_zH_1(\mathbb B,\mathcal D_{T^\ast})\rbrace
\]
and the $K_1$-inner function $W_T$ acts as
\[
W_T(z)=\theta_T(z)|_{\tilde{\mathcal D}}\qquad (z \in \mathbb B).
\]
Hence $\tilde{\mathcal D} \rightarrow W(M),x \mapsto \theta_Tx$, is a unitary operator by Theorem \ref{wandering2}.
\proofend

An elementary exercise shows (cf. Theorem 6.6 in \cite{BEKS}) that, in the setting of the last corollary, the space $\tilde{\mathcal D}$ is also given by
\[
\mathcal D_T\cap({\rm Ker}\ M_{\theta_T})^\perp = \lbrace x \in \mathcal D_T;M^\ast_{\theta_T}M_{\theta_T}x=x\rbrace = \lbrace x \in \mathcal D_T; \| M_{\theta_T}x\|=\| x \|\rbrace.
\]

\vspace{1.5cm}

\centerline{\textbf{\S3 \, $K_m$-Inner functions}}

Let $T\in L(H)^n$ be an $m$-hypercontraction of class $C_{\cdot 0}$ and let $W_T:\mathbb B\rightarrow L(\tilde{\mathcal D},\mathcal D)$,
\[
W_T(z)x = -\tilde{T}x + C\sum^m_{k=1}(1_H-ZT^\ast)^{-k}I^{-1}_mZD_{\tilde{T}}x
\]
be the associated $K_m$-inner function. If we define
\[
B=(\oplus I^{-1}_m)D_{\tilde{T}}:\tilde{\mathcal D}\rightarrow H^n,\quad D=-\tilde{T}:\tilde{\mathcal D}\rightarrow \mathcal D,
\]
then the above representation of $W_T$ becomes
\[
W_T(z)=D+C\sum^m_{k=1}(1_H-ZT^\ast)^{-k}ZB.
\]
The operators $T\in L(H^n,H), B\in L(\tilde{\mathcal D},H^n),C\in L(H,\mathcal D)$ and $D\in L(\tilde{\mathcal D},\mathcal D)$ satisfy the conditions

\begin{align*}
(1)\qquad C^\ast C& = \Delta^{(m)}_T\\ %= \sum_{|\alpha|\leq m}(-1)^{|\alpha|}\binom{m}{|\alpha|}\gamma_\alpha T^\alpha T^{\ast\alpha},\\
(2)\qquad D^\ast C& = -P_{\tilde{\mathcal D}}\tilde{T}^\ast D_{\tilde{T}^\ast} = -P_{\tilde{\mathcal D}}D_{\tilde{T}}\tilde{T}^\ast
=-B^\ast(\oplus I^\ast_mI_m) T^\ast\\
&=-B^\ast (\oplus \sum^{m-1}_{k=0}\Delta^{(k)}_T)T^\ast,\\
(3)\qquad D^\ast D&=P_{\tilde{\mathcal D}}\tilde{T}^\ast \tilde{T}|_{\tilde{\mathcal D}}=1_{\tilde{\mathcal D}}-P_{\tilde{\mathcal D}}(1_{\tilde{H}^n}-\tilde{T}^\ast \tilde{T})|_{\tilde{\mathcal D}}\\
&=1_{\tilde{\mathcal D}}-P_{\tilde{\mathcal D}}D_{\tilde{T}}(\oplus I^\ast_m)^{-1}(\oplus I^\ast_m I_m)(\oplus I_m)^{-1}D_{\tilde{T}}|_{\tilde{\mathcal D}}\\
&=1_{\tilde{\mathcal D}} - B^* (\oplus \sum^{m-1}_{k=0}\Delta^{(k)}_T)B.
\end{align*}

Here $P_{\tilde{\mathcal D}}\in L(\tilde{H}^n,\tilde{\mathcal D})$ denotes the orthogonal projection onto $\tilde{\mathcal D}$. 
Using the definition of the set $\tilde{\mathcal D}\subset \tilde{H}^n$ we obtain the additional condition
\[
(4)\qquad {\rm Im}(\oplus j)B \subset M^\ast_zH_m(\mathbb B,\mathcal D).
\]

Our next aim is to show that the $K_m$-inner functions $W:\mathbb B\rightarrow L(\mathcal E_\ast,\mathcal E)$ between arbitrary Hilbert spaces $\mathcal E_\ast$ and $\mathcal E$ are precisely the 
operator-valued functions on $\mathbb B$ which possess a representation of the form
\[
W(z)=D+C\sum^m_{k=1}(1_H-ZT^\ast)^{-k}ZB,
\]
where $T\in L(H)^n$ is an $m$-hypercontraction of class $C_{\cdot 0}$ on a suitable Hilbert space $H$ and the coefficients of the matrix operator
\[
\biggl(
\begin{array}{c|c}
T^\ast&B\\
\hline
C&D
\end{array}
\biggr)
:H\oplus \mathcal E_\ast \rightarrow H^n \oplus \mathcal E
\]
satisfy the conditions
\begin{align*}
&({\rm KI}1)\ C^\ast C = \Delta^{(m)}_T,\\
&({\rm KI}2)\ D^\ast C = -B^\ast (\oplus \sum^{m-1}_{k=0}\Delta^{(k)}_T)T^\ast,\\
&({\rm KI}3)\ D^\ast D + B^\ast (\oplus \sum^{m-1}_{k=0}\Delta^{(k)}_T)B = 1_{\mathcal E_\ast},\\
&({\rm KI}4)\ {\rm Im}(\oplus j_C)B \subset M^\ast_zH_m(\mathbb B,\mathcal E),
\end{align*}
where $j_C:H\rightarrow H_m(\mathbb B,\mathcal E)$ is the operator defined by
\[
j_Cx=\sum_{\alpha \in \mathbb N^n}\rho_m(\alpha)(CT^{\ast\alpha}x)z^\alpha=C(1_H-ZT^\ast)^{-m}x.
\]
Note that by condition $({\rm KI}1)$ the operator $j_C$ is a well-defined isometry. It is easily seen that $j_C$ intertwines the tuples $T^\ast \in L(H)^n$ and 
$M^\ast_z\in L(H_m(\mathbb B,\mathcal E))^n$ componentwise. The following results extend characterizations of $K_m$-inner functions proved by Olofsson in
\cite{O2} on the unit disc to the case of the unit ball in $\mathbb C^n$.

\begin{thm}
\label{sufficient}
Let $W:\mathbb B \rightarrow L(\mathcal E_\ast,\mathcal E)$ be an operator-valued function between Hilbert spaces $\mathcal E_\ast$ and $\mathcal E$ such that
\[
W(z) = D+C\sum^m_{k=1}(1_H-ZT^\ast)^{-k}ZB,
\]
where $T \in L(H)^n$ is an $m$-hypercontraction of class $C_{\cdot 0}$ and the matrix operator
\[
\biggl(
\begin{array}{c|c}
T^\ast&B\\
\hline
C&D
\end{array}
\biggr)
:H\oplus \mathcal E_\ast \rightarrow H^n \oplus \mathcal E
\]
satisfies the conditions $({\rm KI}1)-({\rm KI}4)$. Then $W$ is a $K_m$-inner function.
\end{thm}

\proof
The intertwining properties of the map $j_C:H\rightarrow H_m(\mathbb B,\mathcal E)$ imply that $M=H_m(\mathbb B,\mathcal E)\ominus ({\rm Im}j_C)\subset H_m(\mathbb B,\mathcal E)$ is a closed $M_z$-invariant
subspace. During the whole proof let $x\in \mathcal E_\ast$ be a fixed vector. According to condition $({\rm KI}4)$ we can choose a function $f \in H_m(\mathbb B,\mathcal E)$ with $(\oplus j_C)Bx=M^\ast_z f$.
Using Lemma \ref{inverse} it follows exactly as in the proof of Theorem \ref{wandering2} that
\[
C\sum^m_{k=1}(1_H-ZT^\ast)^{-k}ZBx = \delta M_z(\oplus j_C)Bx=M_z(M^\ast_zM_z)^{-1}M^\ast_z f.
\]
It follows that
\[
\| C\sum^m_{k=1}(1_H-ZT^\ast)^{-k}ZBx \|^2_{H_m(\mathbb B,\mathcal E)}
\]
\[
= \langle(\oplus j_C)^\ast(M^\ast_zM_z)^{-1}(M^\ast_zM_z)(M^\ast_zM_z)^{-1}M^\ast_z f, Bx\rangle
= \langle(\oplus j_C)^\ast M^\ast_z \delta f,Bx\rangle.
\]
By the remark following Theorem \ref{wandering1} we have
\[
M^\ast_z \delta = (\oplus \sum^{m-1}_{j=0}(-1)^j\binom{m}{j+1}\sum_{|\alpha|=j}\gamma_\alpha M^\alpha_zM^{\ast\alpha}_z)M^\ast_z.
\]
Hence we obtain that
\[
\| C\sum^m_{k=1}(1_H-ZT^\ast)^{-k}ZBx\|^2_{H_m(\mathbb B,\mathcal E)}
\]
\[
= \langle \Bigl( \oplus j^\ast_C (\sum^{m-1}_{j=0}(-1)^j\binom{m}{j+1}\sum_{|\alpha|=j}\gamma_\alpha M^\alpha_zM^{\ast\alpha}_z)j_C \Bigr) Bx,Bx\rangle
\]
\[
= \langle B^\ast \Bigl( \oplus \sum^{m-1}_{j=0}(-1)^j\binom{m}{j+1}\sum_{|\alpha|=j}\gamma_\alpha T^\alpha T^{\ast\alpha} \Bigr) Bx,x\rangle
= \| x \|^2-\|Dx\|^2,
\]
where the last equality follows from condition $({\rm KI}3)$ together with Lemma \ref{defect}. Since $M_zH_m(\mathbb B,\mathcal E)^n=H_m(\mathbb B,\mathcal E)\ominus \mathcal E$, we conclude that the map
\[
\mathcal E_\ast \rightarrow H_m(\mathbb B,\mathcal E),x \mapsto Wx
\]
is a well-defined isometry.\\
Repeating the above calculations we obtain that
\[
M^\ast_z(Wx)=M^\ast_z(M_z(M^\ast_zM_z)^{-1}(\oplus j_C)Bx)=(\oplus j_C)Bx
\]
and hence that
\[
(\oplus j_Cj^\ast_C)M^\ast_z(Wx)=(\oplus j_C)Bx=M^\ast_z(Wx).
\]
Similarly, for $x\in \mathcal E_\ast$ and $f\in H_m(\mathbb B,\mathcal E)$ as above, we find that
\[
j^*_C(Wx)
= C^* Dx+j^*_C(M_z(M^\ast_zM_z)^{-1}M^\ast_z f)
= C^*Dx + T(\oplus j^*_C)M^*_z \delta f
\]
\[
= C^*Dx + T(\oplus j^*_C) \sum^{m-1}_{j=0} (-1)^j \binom{m}{j+1} \sum_{| \alpha | = j} \gamma_{\alpha}M_z^{\alpha} M_z^{* \alpha}M^*_zf
\]
\[
= C^*Dx + T(\oplus \sum^{m-1}_{k=0}\Delta^{(k)}_T)Bx=0,
\]
where the last identity follows from condition $({\rm KI}2)$. Thus we have shown that $W\mathcal E_\ast \subset M\ominus \sum^n_{i=1}z_iM$. Hence
\[
W\mathcal E_\ast \perp z^\alpha(W\mathcal E_\ast)
\]
for all $\alpha \in \mathbb N^n\setminus \lbrace 0 \rbrace$. Thus the proof is complete.
\proofend

%Theorem 13
\begin{thm}
Let $\mathcal E_\ast,\mathcal E$ be Hilbert spaces and let $W:\mathbb B\rightarrow L(\mathcal E_\ast,\mathcal E)$ be a $K_m$-inner function. Then there are a Hilbert space $H$, an $m$-hypercontraction
$T\in L(H)^n$ of class $C_{\cdot 0}$ on $H$ and a matrix operator
\[
\left(
\begin{array}{c|c}
T^\ast&B\\
\hline
C&D
\end{array}
\right)
:H \oplus \mathcal E_\ast \rightarrow H^n \oplus \mathcal E
\]
such that its coefficients satisfy the conditions $({\rm KI}1)-({\rm KI}4)$ and such that
\[
W(z)=D+C\sum^m_{k=1}(1_H-ZT^\ast)^{-k}ZB\quad (z\in \mathbb B).
\]
\end{thm}

\proof
Since $W$ is a $K_m$-inner function, the space $\mathcal W=W(\mathcal E_\ast)\subset H_m(\mathbb B,\mathcal E)$ is a generating wandering subspace for the restriction of $M_z\in L(H_m(\mathbb B,\mathcal E))^n$ to
the closed invariant subspace
\[
\mathcal S=\bigvee_{\alpha \in \mathbb N^n}M^\alpha_z \mathcal W\subset H_m(\mathbb B,\mathcal E).
\]
It is elementary to check and well known that the compression $T=P_HM_z|_H$ of $M_z$ to the co-invariant subspace $H=H_m(\mathbb B,\mathcal E)\ominus \mathcal S$ is an $m$-hypercontraction of class
$C_{\cdot 0}$. Let $\mathcal R\subset H_m(\mathbb B,\mathcal E)$ be the smallest reducing subspace for $M_z\in L(H_m(\mathbb B,\mathcal E))^n$ with $H \subset \mathcal R$. Then
\[
\mathcal R=\bigvee_{\alpha \in \mathbb N^n}z^\alpha (\mathcal R\cap \mathcal E)=H_m(\mathbb B,\mathcal R \cap \mathcal E)
\]
and the inclusion $i:H\hookrightarrow H_m(\mathbb B,\mathcal R \cap \mathcal E)$ is a minimal $m$-dilation for $T$. Let 
$j :H \rightarrow H_m(\mathbb B,\mathcal D)$ be the canonical $m$-dilation of $T$
constructed at the beginning of Section 2. It follows from Section 4 in \cite{K} that there is a unitary operator $U:\mathcal D \rightarrow \mathcal R\cap \mathcal E$ such that 
\[
i = (I \otimes U) \circ j: H = \mathcal S^{\perp} \longrightarrow H_m(\mathbb B, \mathcal R \cap \mathcal E).
\]
Set $\hat{\mathcal E}=\mathcal E\ominus(\mathcal R \cap \mathcal E)$. By definition the space 
\[
H_m(\mathbb B,\hat{\mathcal E})=H_m(\mathbb B,\mathcal E)\ominus H_m(\mathbb B,\mathcal R\cap \mathcal E)=H_m(\mathbb B,\mathcal E)\ominus \mathcal R \subset \mathcal S
\]
is the largest reducing subspace for $M_z\in L(H_m(\mathbb B,\mathcal E))^n$ which is contained in $\mathcal S$. In particular, one obtains the orthogonal decomposition
\[
\mathcal S=H_m(\mathbb B,\hat{\mathcal E})\oplus (\mathcal S\cap H_m(\mathbb B,\hat{\mathcal E})^\perp)=H_m(\mathbb B,\hat{\mathcal E})\oplus (H_m(\mathbb B,\mathcal R\cap \mathcal E)\ominus \mathcal S^\perp). 
\]
Define
\[
M=H_m(\mathbb B,\mathcal D)\ominus {\rm Im}j\in {\rm Lat}(M_z,H_m(\mathbb B,\mathcal D))
\]
and $W(M)=M\ominus(\sum^n_{i=1}z_iM)$. Using the above commutative diagram we find that
\[
I\otimes U:M\rightarrow H_m(\mathbb B,\mathcal R \cap \mathcal E)\ominus \mathcal S^\perp =H_m(\mathbb B,\mathcal R \cap \mathcal E)\cap \mathcal S
\]
defines a unitary operator which intertwines the restrictions of $M_z$ to both sides componentwise. Consequently, we obtain the orthogonal decomposition
\begin{align*}
\mathcal W=\mathcal W(M_z,\mathcal S)&=W(M_z,H_m(\mathbb B,\hat{\mathcal E}))\oplus W(M_z,H_m(\mathbb B,\mathcal R \cap \mathcal E)\cap \mathcal S)\\
&=\hat{\mathcal E}\oplus(I\otimes U)W(M_z,M)=\hat{\mathcal E}\oplus(I\otimes U)W(M).
\end{align*}
Let $W_T:\mathbb B \rightarrow L(\tilde{\mathcal D},\mathcal D)$ be the $K_m$-inner function associated with $T\in L(H)^n$ as in Theorem \ref{wandering2} (see also the beginning of Section 3).
Then there are bounded operators $B \in L(\tilde{\mathcal D},H^n),\ C\in L(H,\mathcal D)$ and $D\in L(\tilde{\mathcal D},\mathcal D)$ such that
\[
W_T(z) = D+C\sum^m_{k=1}(1_H-ZT^\ast)^{-k}ZB\qquad (z \in \mathbb B)
\]
and $W(M) = \lbrace W_T(x);\ x \in \tilde{\mathcal D}\rbrace$. Let us denote by
\[
P_1:\mathcal W\rightarrow \hat{\mathcal E},\ P_2:\mathcal W \rightarrow (I\otimes U)W(M)
\]
the orthogonal projections with respect to the above orthogonal decomposition of $\mathcal W$. The $K_m$-inner functions $W:\mathbb B\rightarrow L(\mathcal E_\ast,\mathcal E)$ and $W_T:
\mathbb B \rightarrow L(\tilde{\mathcal D},\mathcal D)$ induce unitary operators
\[
\mathcal E_\ast \rightarrow \mathcal W,\ x \mapsto Wx,\ \tilde{\mathcal D}\rightarrow W(M),\ x \mapsto W_T(x).
\]
Using these unitary operators we construct bounded linear operators
\[
U_1:\mathcal E_\ast \rightarrow \hat{\mathcal E},\ U_1 x=P_1 Wx
\]
and
\[
U_2:\mathcal E_\ast \rightarrow \tilde{\mathcal D},\ U_2 x=\tilde{x} \mbox{ if } (I\otimes U)W_T\tilde{x}=P_2 Wx.
\]
By construction the column operator $(U_1,U_2):\mathcal E_\ast \rightarrow \hat{\mathcal E}\oplus \tilde{\mathcal D}$ is an isometry such that
\[
W(z)x=U_1x+(I\otimes U)W_T(z)U_2x
\]
\[
=(U_1+UDU_2)x+(UC)\sum^m_{k=1}(1_H-ZT^\ast)^{-k}Z(BU_2)x
\]
for $z \in \mathbb B$ and $x \in \mathcal E_\ast$. To complete the proof it suffices to check that the operators $T\in L(H^n,H),\ \tilde{B}=BU_2\in L(\mathcal E_\ast,H^n),\ \tilde{C}=UC\in L(H,\mathcal E)$ and
$\tilde{D}=U_1+UDU_2 \in L(\mathcal E_\ast,\mathcal E)$ satisfy the conditions $({\rm KI}1)-({\rm KI}4)$. Obviously,
\[
\tilde{C}^\ast\tilde{C}=C^\ast U^\ast UC=C^\ast C=\Delta^{(m)}_T
\]
and
\[
\tilde{D}^\ast \tilde{C}=U^\ast_2D^\ast U^\ast UC=U^\ast_2D^\ast C
\]
\[
=-U^\ast_2B^\ast\left(\oplus \sum^{m-1}_{k=0}\Delta^{(k)}_T\right)T^\ast=-\tilde{B}^\ast\left(\oplus \sum^{m-1}_{k=0}\Delta^{(k)}_T\right)T^\ast.
\]
To check condition $({\rm KI}3)$ note that $\tilde{D}$ acts as the column operator
\[
\tilde{D}=(U_1,UDU_2):\mathcal E_\ast \rightarrow \mathcal E=\tilde{\mathcal E}\oplus (\mathcal R \cap \mathcal E).
\]
Thus we obtain that
\begin{align*}
\tilde{D}^\ast \tilde{D}&=U^\ast_1U_1+U^\ast_2D^\ast U^\ast UDU_2\\
&=U^\ast_1 U_1+U^\ast_2U_2-U^\ast_2B^\ast\left(\oplus \sum^{m-1}_{k=0}\Delta^{(k)}_T\right)BU_2\\
&=I_{\mathcal E_\ast}-\tilde{B}^\ast \left(\sum^{m-1}_{k=0}\Delta^{(k)}_T\right)\tilde{B}.
\end{align*}
Since $j_{\tilde{C}}=Uj_C$, it follows that
\[
(\oplus j_{\tilde{C}})\tilde{B}x=(\oplus U)(\oplus j_C)B(U_2x)\in M^\ast_zH_m(\mathbb B,\mathcal E)
\]
for all $x \in \mathcal E_\ast$. Thus also condition $({\rm KI}4)$ holds. This observation completes the proof.
\proofend

\newpage

\centerline{\textbf{\S4 \, Characteristic functions}}

In the previous section we saw that, for a pure row contraction $T\in L(H)^n$, the associated $K_1$-inner function $W_T:\mathbb B \rightarrow L(\tilde{\mathcal D},\mathcal D)$ is obtained
by restricting its characteristic function
\[
\theta_T:\mathbb B \rightarrow L(\mathcal D_T,\mathcal D_{T^\ast}), \theta_T(z)=-T+D_{T^\ast}(1_H-ZT^\ast)^{-1}ZD_T
\]
to a suitable subspace of $\mathcal D_T$. More precisely,
\[
W_T(z)=\theta_T(z)|_{\tilde{\mathcal D}}\quad (z\in \mathbb B),
\]
where $\tilde{D}\subset D_T$ is the closed subspace given by
\[
\tilde{\mathcal D}=\lbrace y \in \mathcal D_T;\ (\oplus j)D_T y \in M^\ast_zH_1(\mathbb B,\mathcal D_{T^\ast})\rbrace.
\]
In the one-dimensional case $n=1$, we even have the identities $\tilde{\mathcal D}=\mathcal D_T$ and $W_T=\theta_T$. Thus it seens natural to ask whether also for $m> 1$ there is a canonically 
defined characteristic function for each $m$-hypercontraction $T\in L(H)^n$ that extends the function $W_T:\mathbb B\rightarrow L(\tilde{\mathcal D},\mathcal D)$ in a natural way. In the present section we offer
a possible answer to this question.\\

Let $T\in L(H)^n$ be an $m$-hypercontraction of class $C_{\cdot 0}$. We denote by
\[
D_T=(1-T^\ast T)^{1/2}\in L(H^n),\ D_{T^\ast}=(1-TT^\ast)^{1/2}\in L(H)
\]
its  first-order defect operators and by
\[
\mathcal D_T=\overline{D_TH}\subset H^n,\ \mathcal D_{T^\ast}=\overline{D_{T^\ast}H}\subset H
\]
the associated defect spaces. As in Section 2 we write $C = (\Delta^{(m)}_T)^{1/2}\in L(H)$ and $\mathcal D=\overline {CH}\subset H$ for the $m$-th order defect operator and defect space of $T$.

For $k=0,\ldots,m$ the operator
\[
j_k:H\rightarrow H_k(\mathbb B,\mathcal D), j_k(x)=\sum_{\alpha \in \mathbb N^n}\rho_k(\alpha)(CT^{\ast\alpha}x)z^\alpha
\] 
is a well-defined contraction. To check this, it suffices to observe that
\[
\sum_{\alpha \in \mathbb N^n}\frac{\| \rho_k(\alpha)CT^{\ast\alpha}x\|^2}{\rho_k(\alpha)}=\sum_{\alpha \in \mathbb N^n}\rho_k(\alpha)\langle T^\alpha\Delta^{(m)}_TT^{\ast\alpha}x,x\rangle
\]
\[
\leq \sum_{\alpha \in \mathbb N^n}\rho_k(\alpha)\langle T^\alpha \Delta^{(k)}_TT^{\ast\alpha}x,x\rangle=\| x\|^2
\]
for all $x \in H$. Here the above estimate follows from the fact that (Corollary 3 in \cite{MV})
\[
0\leq \Delta^{(m)}_T\leq \Delta^{(m-1)}_T\leq \ldots \leq \Delta^{(1)}_T\leq I.
\]
The inequalities
\[
\langle j^\ast_kj_kx,x\rangle=\sum_{\alpha \in \mathbb N^n}\rho_k(\alpha)\langle T^\alpha\Delta^{(m)}_TT^{\ast\alpha}x,x\rangle \leq \| x\|^2\quad (x \in H)
\]
imply in particular that
\[
j^\ast_kj_k={\rm SOT}-\sum_{\alpha \in \mathbb N^n}\rho_k(\alpha)T^\alpha \Delta^{(m)}_TT^{\ast\alpha}\quad (k=0,\ldots,m).
\]
An elementary calculation yields that
\[
\rho_{k-1}(\alpha)+\sum^n_{\substack{i=1\\ \alpha_i\geq 1}}\rho_k(\alpha-e_i)=\rho_k(\alpha)\quad (k\geq 1,\alpha \in \mathbb N^n).
\]
As a consequence we find that
\[
\sigma_T(j^\ast_kj_k)=\sum^n_{i=1}{\rm SOT}-\sum_{\alpha \in \mathbb N^n}\rho_k(\alpha)T^{\alpha+e_i}\Delta^{(m)}_TT^{\ast \alpha+e_i}
\]
\begin{align*}
&={\rm SOT}-\sum_{\alpha \in \mathbb N^n}\left(\sum^n_{\substack{i=1\\ \alpha_i\geq 1}}\rho_k(\alpha-e_i)\right)T^\alpha \Delta^{(m)}_TT^{\ast\alpha}\\
&=j^\ast_kj_k-j^\ast_{k-1}j_{k-1},
\end{align*}
or equivalently,
\[
(1-\sigma_T)(j^\ast_kj_k)=j^\ast_{k-1}j_{k-1}\quad (k=1,\ldots,m).
\]
Thus we obtain the identities
\[
\Delta^{(k)}_T=(1-\sigma_T)^k(1_H)=(1-\sigma_T)^k(j^\ast_mj_m)=j^\ast_{m-k}j_{m-k}
\]
for $k=0,\ldots,m$.

As before we write $j = j_m : H \rightarrow H_m(\mathbb B,\mathcal D)$. To simplify the notation 
we define $\pi = j_{m-1} : H \rightarrow H_{m-1}(\mathbb B,\mathcal D)$.
Because of
\[
\| \pi(x)\|^2=\| x \|^2-\| T^\ast x\|^2=\| D_{T^\ast} x\|^2\qquad (x \in H)
\]
there is a unique unitary operator $U:\mathcal D_{T^\ast}\rightarrow \mathcal K = \overline{{\rm Im}\pi}$ with $UD_{T^\ast}x=\pi x$ for all $x\in H$. We denote by
\[
\epsilon_z: H_{m-1}(\mathbb B,\mathcal D)\rightarrow \mathcal D, \epsilon_z(f) = f(z)\qquad (z\in \mathbb B)
\]
the point evaluations. Define $M = H_{m-1}(\mathbb B,\mathcal D)\ominus({\rm Im}\pi)$. Then
\[
V=(U,i_M):\mathcal D_{T^\ast}\oplus M\rightarrow H_{m-1}(\mathbb B,\mathcal D),
\]
where $i_M:M\hookrightarrow H_{m-1}(\mathbb B,\mathcal D)$ is the inclusion mapping, is a unitary operator. 
%It is well known that, for a given $m$-hypercontraction $T \in L(H)^n$ of class
%$C_{\cdot 0}$, the mapping
%\[
%j:H\rightarrow H_m(\mathbb B,\mathcal D),j(h)=\sum_{\alpha \in \mathbb N^n}\rho_m(\alpha)(CT^{\ast\alpha}h)z^\alpha=C(1_H-ZT^\ast)^{-m}h
%\]
%is an isometry which intertwines the tuples $T^\ast \in L(H)^n$ and $M^\ast_z\in L(H_m(\mathbb B,\mathcal D))^n$ componentwise. 
Since
\[
\langle h,j^\ast K_m(\cdot,z)x\rangle=\langle j(h)(z),x\rangle=\langle h, (1_H-TZ^\ast)^{-m}Cx\rangle
\]
for all $h \in H,z\in \mathbb B$ and $x \in \mathcal D$, it follows that $j^\ast:H_m(\mathbb B,\mathcal D)\rightarrow H$ is the unique bounded operator with
\[
j^\ast K_m(\cdot,z)x=(1_H-TZ^\ast)^{-m}Cx\qquad (z \in \mathbb B,x \in \mathcal D).
\]
In the case of row contractions, that is $m=1$, we have $\mathcal D=\mathcal D_{T^\ast}$, the map $\pi$ reduces to the mapping $\pi=C=D_{T^\ast}:H\rightarrow 
\mathcal D_{T^\ast}$ and $M=\lbrace 0 \rbrace,U=1_{\mathcal D_{T^\ast}}$. Hence in this case, the operators
\begin{align*}
\Delta_1(z)&=\epsilon_z \circ V:\mathcal D_{T^\ast}\oplus M\rightarrow \mathcal D,\\
\Delta_0&=\epsilon_0 \circ U:\mathcal D_{T^\ast}\rightarrow \mathcal D
\end{align*}
both collapse to the identity operator on $\mathcal D_{T^\ast}$. Thus the only difference between the characteristic functions of row contractions and of $m$-hypercontractions as defined in the 
following theorem is caused by the presence of the space $M$.

\begin{thm}\label{character}
Let $T\in L(H)^n$ be an $m$-hypercontraction of class $C_{\cdot 0}$. Then the operator-valued function $\theta_T:\mathbb B\rightarrow L(\mathcal D_T\oplus M,\mathcal D)$,
\[
\theta_T(z)=-\Delta_1(z)(T \oplus 1_M)+\Delta_0 D_{T^\ast}(1_H-ZT^\ast)^{-m}Z(D_T,0) 
\]
induces a partially isometric multiplier
\[
M_{\theta_T}:H_1(\mathbb B,\mathcal D_T\oplus M)\rightarrow H_m(\mathbb B,\mathcal D)
\]
such that
\[
M_{\theta_T}M^\ast_{\theta_T}+jj^\ast =1_{H_m (\mathbb B,\mathcal D)}.
\]
\end{thm}

\proof
The identity 
\[
\| x \|^2=\| T^\ast x \|^2+\| \pi x \|^2\qquad (x \in H)
\]
obtained above shows that the column operator
\[
H\rightarrow H^n\oplus H_{m-1}(\mathbb B,\mathcal D), x\mapsto (T^\ast x, \pi x)
\]
is an isometry. Define $L=\mathcal D_T\oplus M$ and consider the operators
\begin{align*}
B&=(D_T,0):L=\mathcal D_T\oplus M\rightarrow H^n,\\
D&=-V(T\oplus 1_M):L=\mathcal D_T\oplus M \rightarrow H_{m-1}(\mathbb B,\mathcal D).
\end{align*}
An elementary, but tedious, computation shows that

\[
\biggl(
\begin{array}{c|c}
T^\ast&B\\
\hline
\pi&D
\end{array}
\biggr):H\oplus L\rightarrow H^n \oplus H_{m-1}(\mathbb B,\mathcal D)
\]
defines a unitary matrix operator. Using the identities
\[
\Delta_0D_{T^\ast} = \epsilon_0UD_{T^\ast} = \epsilon_0\pi = C,
\]
one easily obtains the representations
\begin{align*}
\theta_T(z)&=\epsilon_zD+C(1_H-ZT^\ast)^{-m}ZB\\
&=\epsilon_zD + \epsilon_z\pi(1_H-ZT^\ast)^{-1}ZB\\
&=\epsilon_z(D+\pi(1_H-ZT^\ast)^{-1}ZB)
\end{align*}
for all $z\in \mathbb B$. According to Proposition 1.2 in \cite{EP}, the map 
\[
\varphi:\mathbb B \rightarrow L(L,H_{m-1}(\mathbb B,\mathcal D)),\varphi(z)=D+\pi(1_H-ZT^\ast)^{-1}ZB
\]
defines a contractive multiplier from $H_1(\mathbb B,L)$ to $H_1(\mathbb B,H_{m-1}(\mathbb B,\mathcal D))$ such that
\[
K_1(z,w)(1_{H_{m-1}(\mathbb B,\mathcal D)}-\varphi(z)\varphi(w)^\ast)=\pi(1_H-ZT^\ast)^{-1}(1_H-TW^\ast)\pi^\ast
\]
holds for all $z,w \in \mathbb B$. Since
\[
K_m(z,w)1_\mathcal D-\epsilon_z(K_1(z,w)1_{H_{m-1}(\mathbb B,\mathcal D)})\epsilon^\ast_w = (K_m(z,w)-K_1(z,w)K_{m-1}(z,w))1_\mathcal D=0
\]
is a positive definite $L(\mathcal D)$-valued function of $(z,w)\in \mathbb B \times \mathbb B$, the map 
\[
\epsilon: \mathbb B \rightarrow L(H_{m-1}(\mathbb B,\mathcal D), \mathcal D), z \mapsto \epsilon_z
\]
induces a contractive multiplier from $H_1(\mathbb B,H_{m-1}(\mathbb B,\mathcal D))$ to $H_m(\mathbb B,\mathcal D)$. But then also the composition
\[
\mathbb B \rightarrow L(L,\mathcal D), z\mapsto \theta_T(z) = \epsilon_z\varphi(z)
\]
defines a contractive multiplier from $H_1(\mathbb B,L)$ to $H_m(\mathbb B,\mathcal D)$.\\
To complete the proof, note that
\[
\langle jj^\ast K_m(\cdot,w)x,K_m(\cdot,z)y\rangle=\langle 1_H-TW^\ast)^{-m}Cx,(1_h-TZ^\ast)^{-m}Cy\rangle
\]for $z,w \in \mathbb B$ and $x,y \in \mathcal D$. On the other hand,
\begin{align*}
&\langle M_{\theta_T}M^\ast_{\theta_T}K_m(\cdot,w)x,K_m(\cdot,z)y\rangle = \langle \theta_T(w)^\ast x,\theta_T(z)^\ast y\rangle K_1(z,w)\\
&=\langle \epsilon_z \varphi(z)\varphi(w)^\ast \epsilon^\ast_w x,y\rangle K_1(z,w)\\
&=\langle \epsilon_z\epsilon^\ast_w x,y\rangle K_1(z,w)-\langle \epsilon_z(1_{H_{m-1}(\mathbb B,\mathcal D)}-\varphi (z)\varphi(w)^\ast)\epsilon^\ast_w x,y\rangle K_1(z,w)\\
&=K_{m-1}(z,w)K_1(z,w)\langle x,y\rangle-\langle \epsilon_z\pi(1_H-ZT^\ast)^{-1}(1_H-TW^\ast)^{-1}\pi^\ast \epsilon^\ast_w x,y\rangle\\
&=\langle K_m(\cdot,w)x,K_m(\cdot,z)y\rangle - \langle C(1_H-ZT^\ast)^{-m}(1_H-TW^\ast)^{-1}C x,y\rangle
\end{align*}
for all $z,w \in \mathbb B$ and $x,y \in \mathcal D$. Thus
\[
M_{\theta_T}M^\ast_{\theta_T}+jj^\ast=1_{H_m(\mathbb B,\mathcal D)}
\]
and the proof is complete.
\proofend

Since the matrix operator
\[
\left(
\begin{array}{c|c}
T^\ast&B\\
\hline
\pi&D
\end{array}
\right):H\oplus L\rightarrow H^n\oplus H_{m-1}(\mathbb B,\mathcal D)
\]
is unitary, its second column defines a unitary operator $\rho:L\rightarrow K$ onto the orthogonal complement $K$ of the image of its first column in $H^n\oplus H_{m-1}(\mathbb B,\mathcal D)$.
Using the definitions of $B$ and $D$ we find that
\[
\rho(x,f)=(D_T x,-UTx-f)
\]
for $(x,f)\in L=\mathcal D_T\oplus (H_{m-1}(\mathbb B,\mathcal D)\ominus {\rm Im}\ j_{m-1})$. Let us denote by
\[
\mathcal B:K\rightarrow H^n,(x,f)\mapsto x,\mathcal D:K\rightarrow H_{m-1}(\mathbb B,\mathcal D),(x,f)\mapsto f
\]
the projections of $K$ onto its first and second component. Then
\[
\varphi:\mathbb B \rightarrow L(K,\mathcal D),\varphi (z)=\epsilon_z \circ \mathcal D+C(1_H-ZT^\ast)^{-m}Z\mathcal D
\]
defines an analytic operator-valued function such that
\[
\varphi(z)\circ \rho=\theta_T(z)\quad (z\in \mathbb B).
\]
Hence also $M_\varphi:H_1(\mathbb B,K)\rightarrow H_m(\mathbb B,\mathcal D)$ is a partially isometric multiplier and
\[
M_\varphi M^\ast_\varphi + jj^\ast=1_{H_m(\mathbb B,\mathcal D)}.
\]

Using the above alternative characterization of $\theta_T$, one can show that $\theta_T$ is a purely contractive multiplier (cf. Section V.2 in \cite{Sz-N} and Section 4 in \cite{BES2}).

\begin{lem}
Let $T \in L(H)^n$ be an $m$-hypercontraction of class $C_{\cdot 0}$. Then its characteristic function $\theta_T:\mathbb B\rightarrow L(\mathcal D_T\oplus M,\mathcal D)$ satisfies
\[
\| \theta_T(0)x\| < \| x \|
\]
for all non-zero vectors $x \in \mathcal D_T \oplus M$.
\end{lem}

\proof
It suffices to prove the corresponding result for $\varphi$. Note that the operator $\varphi (0):K\rightarrow \mathcal D$ acts as
\[
\varphi(0)(x,f)=f(0).
\]
So, if $(x,f)\in K$ is a vector with $\| \varphi (0)(x,f)\| \geq \|(x,f)\|$, then
\[
\| (x,f)\|=\| f(0)\| \leq \| f \|_{H_{m-1}(\mathbb B,\mathcal D)}.
\]
Hence $x =0$ and $\| f(0)\|=\| f \|_{H_{m-1}(\mathbb B,\mathcal D)}$. But then $f\equiv f(0)$ is a constant function. Since $(x,f)$ is orthogonal to the image of the column operator
\[
H \rightarrow H^n \oplus H_{m-1}(\mathbb B,\mathcal D), h\mapsto (T^\ast h,j_{m-1}h),
\]
we conclude that $f\in H_{m-1}(\mathbb B,\mathcal D)\ominus {\rm Im}\ j_{m-1}={\rm Ker}\ j^\ast_{m-1}$. Since
\[
\langle h,j^\ast_{m-1}z\rangle=\langle j_{m-1}h,z\rangle=\langle (j_{m-1}h)(0),z\rangle=\langle Ch,z\rangle=\langle h,Cz\rangle
\]
for all $h \in H$ and $z\in \mathcal D$, it follows that $j^\ast_{m-1}z=Cz$ for $z \in \mathcal D$. Thus we find that also
\[
f\equiv f(0)\in \mathcal D \cap {\rm Ker}\ C=(\overline{{\rm Im}\ C})\cap ({\rm Ker}\ C)=\lbrace 0 \rbrace.
\]
This observation completes the proof.
\proofend

J\"org Eschmeier\\
Fachrichtung Mathematik\\
Universit\"at des Saarlandes\\
Postfach 15 11 50\\
D-66041 Saarbr\"ucken,
Germany\\
e-mail: eschmei@math.uni-sb.de\\

\end{document}